\documentclass[11pt]{article}
\usepackage{amsmath}
\usepackage{amssymb}
\usepackage{amsthm}
\usepackage{amsfonts}
\usepackage{graphicx}
\usepackage{pdfpages}
\usepackage{xcolor,colortbl}
\usepackage{multicol}
\usepackage{url}
\usepackage{subfigure}

\definecolor{Gray}{gray}{0.92}
\definecolor{LightCyan}{rgb}{0.88,1,1}

\newcolumntype{a}{>{\columncolor{Gray}}c}
\newcolumntype{b}{>{\columncolor{white}}c}

\usepackage{bbm,epsfig,graphics,epic,color,rotating,color}
    \usepackage[left=2.5cm,right=2.5cm,top=2.3cm,bottom=2.5cm]{geometry}

\numberwithin{equation}{section}
\def\E{{\mathbb E}}

\def\P{{\mathbb P}}

\newtheorem*{theorem*}{Theorem}
\newtheorem{theorem}{Theorem}[section]
\newtheorem{lemma}[theorem]{Lemma}

\newtheorem{definition}[theorem]{Definition}

\newtheorem{remark}[theorem]{Remark}

\theoremstyle{definition}

\title{Reinforced random walks under memory lapses}
\author{Manuel Gonz\'alez-Navarrete and Ranghely Hern\'andez}
\date{}

\begin{document}

\maketitle

\vspace{2pt}

\begin{abstract}

We introduce a one-dimensional random walk, which at each step performs a reinforced dynamics with probability $\theta$ and with probability $1 - \theta$, the random walk performs a step independent of the past. We analyse its asymptotic behaviour, showing a law of large numbers and characterizing the diffusive and superdiffusive regions. We prove central limit theorems and law of iterated logarithm based on the martingale approach.
\end{abstract}

\section{Introduction}

In the area of mathematical physics, the study of random dynamics that combine simplicity in their definition and complexity in their behaviors has drawn vital relevance. In this sense, the so-called elephant random walk (ERW) has generated a lot of interest in recent years. The ERW was proposed by Schutz and Trimper \cite{schutz2004elephants} and can be formulated as follows. Let the sequence $\{X_1, X_2, \ldots\}$ where $X_i \in \{-1,+1\}$, for all $i\ge 1$. The position of the elephant at instant $n$ is given by $S_{n} = \sum_{i = 1}^{n} X_{i},$ and $S_0=0$. It is assumed that the elephant completely remembers its past and performs the $(n+1)$-step by choosing $t\in \{1,....n\}$ uniformly at random, then the elephant decides to repeat the step $X_t$ with probability $p$ or to do the opposite with probability $1-p$.

The main characteristic of ERW is the existence of a phase transition from a diffusive to a super-diffusive behavior (\cite{baur2016elephant,coletti2017central}). This refers to the fact that the second moment of the position is scaled linearly (diffusive) or as a power law (super-diffusive). This fact has motivated the development of new models and results that generalize the ERW, for instance, Bercu and Laulin \cite{bercu2019multi,BL2}, Bertoin \cite{bertoin2018noise,bertoin2020universality,Bert2}, González-Navarrete and Lambert \cite{gonzalez2020multidimensional,gonzalez2019diffusion} and Gut and Stadtmüller \cite{gut2019elephant,gut2020zeros}.

This work is motivated by the recent formulation in \cite{gut2019elephant}, which defines an ERW with delays. In particular, we propose the addition of memory lapses to the dynamics of the process, this kind of evolution was introduced by \cite{gonzalez2018non} and also used in \cite{Baur,Busi}. In this sense, at each step, with probability $\theta$ the elephant performs its movement by following the dynamics in \cite{gut2019elephant}, otherwise the elephant performs a step independent of its history.

We aim to characterize the asymptotic behaviors of the proposed model. In particular, it seeks to analyse the influence of the parameters $p \in [0,1]$ and $\theta \in [0,1)$. Specifically, we prove a law of large numbers and characterize the diffusive and superdiffusive regions, showing convergences in distribution and laws of iterated logarithm. The main technique we used is the martingale theory, based on the works \cite{bercu2017martingale,duflo2013random,hall2014martingale,heyde1977central,JJQ}.

The rest of the paper is organized as follows. In Section \ref{cpt_mod} we include the mathematical formulation of the model and state the main results. Section \ref{sec:proof} includes some auxiliary results and the proofs of the theorems.

\section{The model and main results}
\label{cpt_mod}

Gut and Stadtm\"uller \cite{gut2019elephant} proposed an ERW with memory defined as follows. Let the first step given by
\begin{equation}
\label{paso1}
X_1=\left\{
\begin{array}{ll}
+1   &\  \mbox{, with probability $p$,} \\
-1  &\ \mbox{, with probability $q$,}\\
0  &\ \mbox{, with probability $r$.}\\
\end{array}
\right.
\end{equation}
where $p+q+r=1$. The next steps are performed by the rule
\begin{equation}
\label{original}
X_{n+1} = \left\{
\begin{array}{lll}
X_{t}      &, \  \mbox{with probability $p$,} \\
-X_{t} &, \ \mbox{with probability $q$,} \\
0 &, \ \mbox{with probability $r$,}
\end{array}
\right.
\end{equation}
where $t$ is uniformly distributed on $\{1,2,\ldots\,n\}$.
The sequence $\{X_i\}$ generates a one-dimensional random walk $S_n$ given by:
\begin{equation}
    \label{Sn}
    \begin{array}{ccccc}
       S_0:=0 & \text{ and } &S_n=\displaystyle\sum_{i=1}^{n} X_i & \forall ~ n=1,2,\ldots
    \end{array}
\end{equation}
In this sense, \cite{gut2019elephant} proved for $0 \leq p - q < 1/2$ that
    \begin{equation}
    \label{GS1}
        \frac{S_n}{\sqrt{n}}\overset{d}{\rightarrow} (p+q)N\left(0,\frac{p+q}{1-2(p-q)}\right) + r \delta_0(x), \ \text{ as } n \to \infty
    \end{equation}
where $\overset{d}{\rightarrow}$ means convergence in distribution. Moreover, if $p - q =1/2$, then
    \begin{equation}
    \label{GS2}
        \frac{S_n}{\sqrt{n\cdot \log n}}\overset{d}{\rightarrow} (p+q)N\left(0,p+q\right) + r \delta_0(x), \ \text{ as }  n \to \infty.
    \end{equation}
we use the standard $\delta_0(x)$ to denote the distribution function with a jump of height one at $0$. The results in \eqref{GS1} and \eqref{GS2} are given by the fact that, if $X_1 = 0$, then $S_n = 0$, $\forall ~ n \ge 1$, which occurs with probability $r$.
As usual, in the region $1/2<p - q <1$, there exists a random variable with positive variance, such that $\lim\limits_{n\to\infty}\frac{S_n}{n^{p-q}}=W$ almost surely as $n \to \infty$.

Based on \eqref{original} we propose a model with random tendency as described in \cite{gonzalez2020multidimensional}. At each step $n\geq 1$, we draw a Bernoulli random variable $Y_n$ with parameter $\theta$, such that for $Y_n=0$ the random walk performs a step independent of its past, otherwise if $Y_n=1$ the elephant follows the dynamics proposed in \eqref{original}. The conditional probabilities can be formulated as follows. Let $\mathcal{F}_n$ the $\sigma$-algebra generated by $\{X_1, X_2,\ldots,X_n\}$, for all $n=1,2,\ldots$
\begin{equation}
\label{PT1}
\begin{array}{c}
\P(X_{n+1}=+1\vert \mathcal{F}_n) = \left( \dfrac{n_+}{n}\cdot p+\dfrac{n_-}{n}\cdot q\right)\cdot \theta+(1-\theta)\cdot p \\[0.5cm]
\P(X_{n+1}=-1\vert \mathcal{F}_n) = \left( \dfrac{n_-}{n}\cdot p+\dfrac{n_+}{n}\cdot q\right)\cdot \theta+(1-\theta)\cdot q \\[0.5cm]
\P(X_{n+1}=0\vert \mathcal{F}_n) = \left( \dfrac{n_0}{n}\cdot (p+q)+r\right)\cdot\theta+(1-\theta)\cdot r = \dfrac{n_0}{n}\cdot (p+q)\cdot \theta+r 
\end{array}
\end{equation}
where $\theta\in [0,1)$, $n_+ = \#\{ i \in \{1, 2, \ldots, n\}: X_i = +1\}$, is the number of $+1$ steps until the $n$-step, and $n_-$ and $n_0$ represent the number of steps in the directions $-1$ and $0$ until the $n$-step, respectively.

\subsection{Main results}

We use the following notation
\begin{equation}
\label{notation}
\alpha=(p-q)\cdot\theta, \ \omega=(p-q)(1-\theta), \ \tau=(1-\theta)(p+q), \ \gamma=(p+q)\theta \text{ and } \phi=\frac{\tau}{1-\gamma}-\left(\frac{\omega}{1-\alpha}\right)^2
\end{equation}
and let denote the sequence of constants
\begin{equation}
\label{an}
\begin{array}{cccc}
a_1:=1& \text{ and } a_n:=\displaystyle\prod_{k=1}^{n-1}\left(1+\dfrac{\alpha}{k}\right)=\dfrac{\Gamma(n+\alpha)}{\Gamma(n)\cdot\Gamma(\alpha+1)}&\text{ for }& n=2,\ldots
\end{array}
\end{equation}

The main results are

\begin{theorem}
\label{Teo1}
Let the RW given by \eqref{PT1}, for all $\alpha\in [0,1)$
\begin{equation}
\label{lgn}
\begin{array}{lll}
   \lim\limits_{n\to\infty}\frac{S_n-\E(S_n)}{n}=0  & \text{ a.s}
\end{array}
\end{equation}
and
\begin{equation}
\label{limSn}
\begin{array}{ll}
  \lim\limits_{n\to\infty}\frac{S_n}{n}=\frac{\omega}{1-\alpha}   &  a.s.
\end{array}
\end{equation}
\end{theorem}

%Donde $\omega=(1-\theta)(p-q)$ y $\alpha=(p-q)\cdot\theta$.

\begin{theorem}
\label{T_conv}
Let the RW given by \eqref{PT1}:
\begin{itemize}
    \item [i)] If $0 \leq \alpha < 1/2$, then
    \begin{equation}
    \label{TCLdif}
        \frac{S_n-\E[S_n]}{\sqrt{\frac{\phi}{1-2\alpha}n}}\overset{d}{\rightarrow}N(0,1), \text{ as } \ n \to \infty
    \end{equation}
    and
    \begin{equation}
    \label{LIL1}
        \limsup_{n \to \infty} \pm \frac{S_n-\E[S_n]}{\sqrt{2 \frac{\phi}{1-2\alpha} n \log \log \left(\frac{\phi \Gamma(\alpha + 1)^2 }{1-2\alpha}n^{1-2\alpha}\right)}} = 1 \ \text{ a.s }
    \end{equation}
%    Donde $\phi=\frac{\tau}{1-\gamma}-\left(\frac{\omega}{1-\alpha}\right)^2$, $\tau=(1-\theta)(p+q)$ y $\gamma=(p+q)\theta$

    \item[ii)] If $\alpha=1/2$, then
    \begin{equation}
        \label{TCLcrit}
        \frac{S_n-\E[S_n]}{\sqrt{\phi\cdot n \cdot \log n}}\overset{d}{\rightarrow}N(0,1), \text{ as }  \ n \to \infty
    \end{equation}
      and
    \begin{equation}
    \label{LIL2}
        \limsup_{n \to \infty} \pm \frac{S_n-\E[S_n]}{\sqrt{2 \phi \cdot  n  \cdot \log n \cdot \log \log (\phi\Gamma(3/2) \log n)}} = 1 \ \text{ a.s }
    \end{equation}
  %   Donde $\phi=\frac{\tau}{1-\gamma}-\left(\frac{\omega}{1-\alpha}\right)^2$, $\tau=(1-\theta)(p+q)$ y $\gamma=(p+q)\theta$
  
\end{itemize}
\end{theorem}

\begin{theorem}
If $1/2<\alpha<1$, then there exists a random variable $W$ with positive variance such that
    \begin{equation}
        \label{W}
        \lim\limits_{n\to\infty}\frac{S_n-\E[S_n]}{a_n}=W \ \text{a.s.} \ \text{and in } L^2,
    \end{equation}
    Moreover,

\begin{equation}
    \label{T_cnv_w}
    \frac{S_n-\E[S_n]-W\cdot a_n}{\sqrt{\frac{\phi}{2\alpha-1}n}}\overset{d}{\rightarrow}N(0,1), \text{ as }  \ n \to \infty
\end{equation}
   and
    \begin{equation}
    \label{LIL3}
        \limsup_{n \to \infty} \pm \frac{S_n-\E[S_n]-W\cdot a_n}{\sqrt{2 \frac{\phi}{2\alpha-1}n \log \left| \log \left(\frac{\phi\Gamma(\alpha+1)^2}{2\alpha-1}n^{1-2\alpha}\right)\right|}} = 1 \ \text{ a.s }
    \end{equation}
%Donde  $\phi=\frac{\tau}{1-\gamma}-\left(\frac{\omega}{1-\alpha}\right)^2$ y $\alpha=(p-q)\cdot\theta$.
\end{theorem}

We remark that the distribution and higher moments of $W$ are of interest. In the case of original ERW we refer \cite{bercu2017martingale,hyper} for details.

\section{Proofs}
\label{sec:proof}

First, we need to introduce some definitions and auxiliary results

\begin{definition}
Let define the sequence
\begin{equation}
\label{Mart}
\begin{array}{ccc}
M_n=\dfrac{S_n-\E(S_n)}{a_n}&\text{ for }& n=1,2,\ldots
\end{array}
\end{equation}
where $a_n$ as defined in \eqref{an}.
\end{definition}
We will prove that $M_n$ is a martingale (see Lemma \ref{l_mart}) and the asymptotic behaviour of the martingale is strictly related with the sum:
\begin{equation}
    \label{vn}
    v_n=\displaystyle\sum_{k=1}^{n}\frac{1}{a_k^2}
\end{equation}
Now, note that by Stirling formula for the gamma function
\begin{equation}
\label{conv_an}
\begin{array}{cccc}
a_n \thicksim \frac{n^\alpha}{\Gamma(\alpha+1)}&\text{ as }& n\rightarrow \infty
\end{array}
\end{equation}
Therefore, in the \textbf{diffusive region}, where $0\leq\alpha<1/2$, we will prove in Lemma \ref{lan2} that $v_n$ diverges. In particular, by \eqref{an} we have:
\begin{equation}
    \label{conv_dif_vn}
        v_n =\displaystyle\sum_{k=1}^{n}\left(\frac{\Gamma(k)\Gamma(\alpha+1)}{\Gamma(k+\alpha)}\right)^2  \sim (\Gamma(\alpha+1))^2\displaystyle\sum_{k=1}^{n}\frac{1}{k^{2\alpha}}
\end{equation}
as $n\rightarrow\infty$. Then, by the p-series we know that $v_n  \sim(\Gamma(\alpha+1))^2\cdot\frac{n^{1-2\alpha}}{1-2\alpha}$. Then,
  \begin{equation}
  \label{A8}
  \lim\limits_{n\to\infty}\frac{v_n}{n^{1-2\alpha}}=\frac{(\Gamma(\alpha+1))^2}{1-2\alpha}
  \end{equation}
In the \textbf{critical region}, where $\alpha=1/2$, we have $ v_n  \sim (\Gamma(3/2))^2\log n$ and
\begin{equation}
    \label{conv_crit_vn}
    \begin{array}{ll}
        v_n & \sim (\Gamma(3/2))^2\displaystyle\sum_{k=1}^{n}\frac{1}{k}
    \end{array}
\end{equation}
Then $v_n$ diverges with velocity $\log n$ (see Lemma \ref{lan2}), by the harmonic series, we obtain
\begin{equation}
    \label{A9}
    \lim\limits_{n\to\infty}\frac{v_n}{\log n}=\frac{\pi}{4}
\end{equation}
Finally, in the \textbf{superdiffusive region}, if $1/2<\alpha\leq1$, $v_n$ converges to the finite value
\begin{equation}
    \label{A10}
    \begin{array}{ll}
      \displaystyle\lim_{n\to\infty}v_n &=\displaystyle\sum_{k=0}^{\infty}\left(\frac{\Gamma(k+1)\Gamma(\alpha+1)}{\Gamma(k+\alpha+1)}\right)^2=\displaystyle\sum_{k=0}^{\infty}\frac{(1)_k(1)_k(1)_k}{(\alpha+1)_k(\alpha+1)_k k!}\\[0.3cm]
      &={}_3 F_2 \left(
    \begin{array}{c|c}
        1,1,1& \\  
        (\alpha+1),(\alpha+1) &
    \end{array} 1 \right)
    \end{array}
\end{equation}
where for $a\in\mathbb{R}$, $(a)_k=a(a+1)\cdots(a+k)$ for $k\geq0$ and $(a)_0=1$ represent the Pochhammer symbol and ${}_3 F_2$ the hypergeometric generalized function.

We also will use the following definitions.

%%%%%%%%%%%%%%%%%%%%%%%%%%%%%%%%%%%%%%%%%%%%%%%%%%%%%%%%

\begin{definition}
The predictable quadratic variation associated to $(M_n)$ is given by $\langle M\rangle_0=0$ and, for all $n\geq 1$
\begin{equation}
    \label{Var_M}
    \langle M\rangle_n=\displaystyle\sum_{k=1}^{n}\E[\Delta M^2_k\vert \mathcal{F}_{k-1}].
\end{equation}
\end{definition}

\begin{definition}
Let denote $Z_n=\displaystyle\sum_{i=1}^{n}(X_i)^2$, $\forall ~ n=1,2,\ldots$, and $\{M^*_n\}=\frac{Z_n-\E[Z_n]}{b_n}$ where $b_n$ is given by $b_1=1$ and $b_n:=\displaystyle\prod_{k=1}^{n-1}\left(1+\dfrac{\gamma}{k}\right)$ for all $n\ge 2$, with $\gamma$ as in \eqref{notation}.
\end{definition}

Now, we prove some lemmas useful in the proofs of the Theorems.

\begin{lemma}
\label{l_mart}
The sequences $\{M_n\}$ and $\{M^*_n\}$ are square integrable martingales.
\end{lemma}

\textbf{Proof}
First we prove that $\{M_n\}$ is a martingale, in fact from \eqref{PT1}
\begin{equation}
\label{XndadoF}
\begin{array}{lll}
\E(X_{n+1}\vert S_n)&=& 1\cdot \P(X_{n+1}=1\vert \mathcal{F}_n)-1\cdot \P(X_{n+1}=-1\vert \mathcal{F}_n)+0\cdot\P(X_{n+1}=0\vert \mathcal{F}_n)\\[0.5cm]
    
       &=&\dfrac{(p-q)\cdot \theta}{n}\cdot S_n +(1-\theta)\cdot(p-q)
\end{array}
\end{equation}
and
\begin{equation}
    \label{Esp}
    \E(S_{n+1}\vert S_n) = \lambda_n S_n+\omega
\end{equation}
where $\lambda_n=1+\dfrac{\alpha}{n}$ with $\alpha=(p-q)\cdot\theta$ and $\omega=(1-\theta)\cdot (p-q)$.
\begin{equation}
    \label{mart}
    \begin{array}{ll}
      \E(M_{n+1}\vert\mathcal{F}_n) & = \E\left[\frac{S_n+X_{n+1}-\E(S_n+X_{n+1})}{a_{n+1}}\vert\mathcal{F}_n\right]\\[0.3cm]
          &=\frac{(n+\alpha)(S_n-\E(S_n))}{n\cdot\left(1+\frac{\alpha}{n}\right)a_{n}} = M_n
          \end{array}
\end{equation}
Analogously, note that $ \E(M^*_{n+1}\vert\mathcal{F}_n)$ equals
\begin{equation}
    \label{martZ}
    \begin{array}{ll}
      \E\left[\frac{Z_{n+1}-\E(Z_{n+1})}{b_{n+1}}\vert\mathcal{F}_n\right]
        & = \E\left[\frac{Z_n+X^2_{n+1}-\E(Z_n+X^2_{n+1})}{b_{n+1}}\vert\mathcal{F}_n\right]\\[0.3cm]
          & =\frac{nZ_n+\gamma Z_n-n\E(Z_{n})-\gamma \E(Z_n)}{nb_{n+1}} =M^*_n
        \end{array}
\end{equation}
In addition, for any $n \geq 1$, $X_n$ are random variables with values on $\{+ 1,0,-1\}$. Then, $|S_n|, |Z_n|\leq n$, which implies that $(M_n)$ are $(M^*_n)$ are locally square integrable martingales.\\

%%%%%%%%%%%%%%%%%%%%%%%%%%%%%%%%%%%%%%%%%%%%%%%%%%%%%%%%%%%%%%%%%%

\begin{lemma}
The martingale $(M_n)$ can be written in the additive form
\begin{equation}
\label{adi_M}
    M_n=\displaystyle\sum_{k=1}^{n-1} \Delta M_{k}= \displaystyle\sum_{k=1}^{n-1} (M_{k}-M_{k-1})=\displaystyle\sum_{k=1}^{n-1}\frac{X_k-\E(X_k\vert \mathcal{F}_{k-1})}{a_k}
\end{equation}
\end{lemma}

\textbf{Proof}
First we calculate the expectation
\begin{equation}
    \label{dif}
    \begin{array}{c}
       \E(S_{n}-\lambda_{n-1}S_{n-1})   =\E\left[S_{n-1}+X_n-\left(1+\frac{\alpha}{n-1}\right)S_{n-1}\right]\\[0.4cm]
          =\E\left[\E\left[X_n\vert S_{n-1}\right]-\frac{\alpha}{n-1}S_{n-1}\right] =\omega        \end{array}
\end{equation}
Now, the martingale difference $\Delta M_{n}=M_{n}-M_{n-1}$, using \eqref{XndadoF} and \eqref{dif}:

\begin{equation*}
    \begin{array}{ll}
   \Delta M_n &=\frac{S_{n}-\E(S_{n})}{a_{n}}- \frac{S_{n-1}-\E(S_{n-1})}{a_{n-1}}
     =\frac{S_{n}-\lambda_{n-1}S_{n-1}-(\E[S_{n}-\lambda_{n-1} S_{n-1}])}{a_{n}}\\[0.4cm]
     &=\frac{X_n-\E(X_n\vert \mathcal{F}_{n-1})}{a_{n}}
     \end{array}
\end{equation*}

\begin{remark}
\label{cota_dk}
Given that $X_n$ assumes values on $\{-1,0,1\}$ it is easy to see that
\begin{equation}
    \label{AcMn}
    \frac{\vert \Delta M_n\vert}{a_n}\leq \frac{2}{a_n}
\end{equation}
\end{remark}

%%%%%%%%%%%%%%%%%%%%%%%%%%%%%%%%%%%%%%%%%%%%%%%%%%%%%

%%%%%%%%%%%%%%%%%%%%%%%%%%%%%%%%%%%%%%%%%%%%%%%%%%%%%%%%%%%%
\begin{lemma}
For $n=1,2,\ldots,$
\begin{equation}
\label{lEsn}
    \E[S_n]=\beta a_n+\omega\cdot a_n \displaystyle\sum_{l=1}^{n-1}\frac{1}{a_{l+1}}
\end{equation}
and
\begin{equation*}
\label{lzn}
    \E[Z_n]=\psi b_n+\tau\cdot b_n \displaystyle\sum_{l=1}^{n-1}\frac{1}{b_{l+1}}
\end{equation*}
\end{lemma}
\textbf{Proof}\\
Note that:
 $x_1=\E[S_1]=\E[X_1]=p-q=\beta$, and $x_{n+1}=E[S_{n+1}]=\E[\E[S_{n+1}\vert S_n]]=\lambda_n \E(S_n)+\omega$ by \eqref{Esp}. Now, using $x_1=f_0=\beta$, $f_n=\omega, g_n=\lambda_n$ and $x_n=\E(S_n)$ for $n=1,2,\ldots,$ we have by Lemma A1 from \cite{kubota2019gaussian} that:
\begin{equation*}
    \E[S_n]=\beta\cdot \displaystyle\prod_{k=1}^{n-1}\lambda_n+\omega\displaystyle\sum_{l=1}^{n-1}\displaystyle\prod_{k=l+1}^{n-1}\lambda_k \\[0.5cm] 
          = \beta a_n + \omega\cdot a_n\displaystyle \sum_{l=1}^{n-1}\frac{1}{a_{l+1}}.
\end{equation*}
In a similar form we obtain $ \E[Z_n]$.

\begin{lemma}
\label{lemma_Sum_an}
Let $a_n$ defined in \eqref{an}, then $\displaystyle\sum_{l=1}^{n-1}\frac{1}{a_{l+1}}\sim\frac{\Gamma(\alpha+1)n^{1-\alpha}}{(1-\alpha)}$.
\end{lemma}

\textbf{Proof}
Note that
\begin{equation*}
       \displaystyle\sum_{l=1}^{n-1}\frac{1}{a_{l+1}}= \displaystyle\sum_{l=2}^{n}\frac{1}{a_{l}}
    =\displaystyle\sum_{l=1}^{n}\frac{\Gamma(l)\Gamma(\alpha+1)}{\Gamma(l+\alpha)}-\frac{1}{a_1}
    =\Gamma(\alpha+1)\displaystyle\sum_{l=1}^{n}\frac{\Gamma(l)}{\Gamma(l+\alpha)}-1
\end{equation*}
By using Lemma B.1 in \cite{bercu2017martingale}, that is $\sum_{l=1}^{n}\frac{\Gamma(l)}{\Gamma(l+\alpha)}= \frac{\Gamma(n+1)}{(\alpha-1)\Gamma(n+\alpha)}\cdot\left(\frac{\Gamma(n+\alpha)\Gamma(1)}{\Gamma(n+1)\Gamma(\alpha)}-1\right)$, the equation \eqref{an} and taking $a=0$, $b=\alpha$, with $\alpha\neq1$, it holds:
\begin{equation*}
       \displaystyle\sum_{l=1}^{n-1}\frac{1}{a_{l+1}}
        =\frac{n}{(1-\alpha)a_n}+\frac{1}{\alpha-1}\\[0.5cm] 
       \sim\frac{\Gamma(\alpha+1)n^{1-\alpha}}{(1-\alpha)} 
\end{equation*}
%%%%%%%%%%%%%%%%%%%%%%%%%%%%%%%%%%%%%%%%%%%%%%%%%%
\begin{lemma}
\label{lan2}
The series $\sum\frac{1}{a_n^2}$ converges, if and only if, $\alpha>\frac{1}{2}$
\end{lemma}

\textbf{Proof}
Let $C_n=\frac{1}{a_n^2}$ by the Raabe criteria, we need to analyse $\lim\limits_{n\to\infty} n\left(1-\frac{C_{n+1}}{C_n}\right)$. Note that,
 \begin{equation*}
     n\left(1-\frac{C_{n+1}}{C_n}\right)=n\left(1-\frac{1}{\lambda_n^2}\right)=n\left(\frac{\frac{2\alpha}{n}+\frac{\alpha^2}{n^2}}{1+\frac{2\alpha}{n}+\frac{\alpha^2}{n^2}}\right)=2\alpha
\end{equation*}
The series converges if $\alpha>1/2$, and diverges if $\alpha<1/2$. In the case $\alpha=1/2$, from \eqref{an} it follows from \eqref{conv_an} that $a_n\thicksim \frac{n^{1/2}}{\sqrt{\pi/4}}$. Since $\frac{1}{a_n^2}=\frac{\pi/4}{n}\geq\frac{1}{2n}$, and the series $\sum\frac{1}{a_n^2}$ diverges.

Finally, we state the following lemma without its proof. The details can be obtained in Lemma \ref{Kubotalem} in \cite{kubota2019gaussian}, and the proof is based on Theorem 1 (b) from \cite{heyde1977central}.

\begin{lemma}
\label{Kubotalem}
Suppose that $\{M_n\}$ is a square-integrable martingale with mean $0$. Let $\Delta M_k = M_k - M_{k-1}$, for $k=1, 2, \ldots$, where $M_0=0$. If
\begin{equation*}
    \displaystyle\sum_{k=1}^{\infty} \E[(\Delta M_k)^2] < +\infty
\end{equation*}
holds in addition, then we have the following: let $r_n^2= \displaystyle\sum_{k=n}^{\infty} \E[(\Delta M_k)^2] $
\begin{itemize}
    \item [(i)] The limit $M_{\infty}:= \sum_{k=1}^{\infty} \Delta M_k $ exists almost surely and $M_n \overset{L^2}{\rightarrow} M_{\infty}$
    \item [(ii)] Assume that
    \begin{itemize}
        \item [a)] $\displaystyle\frac{1}{r_n^2}\sum_{k=n}^{\infty} (\Delta M_k)^2 \to 1$ as $n \to \infty$ in probability, and
        \item [b)] $\displaystyle\lim_{n\to\infty}\frac{1}{r_n^2}\E\left[ \sup_{k \ge n} \ (\Delta M_k)^2 \right] =0$.
        Then, we have
        \begin{equation*}
            \displaystyle\frac{M_{\infty}-M_n}{r_{n+1}}=\frac{\sum_{k=n+1}^{\infty}\Delta M_k}{r_{n+1}} \overset{d}{\rightarrow}N(0,1).
        \end{equation*}
    \end{itemize}
     \item [(iii)] Assume that the following three conditions hold 
       \begin{itemize}
        \item [a')] $\displaystyle\frac{1}{r_n^2}\sum_{k=n}^{\infty} (\Delta M_k)^2 \to 1$ as $n \to \infty$ a.s
        \item [c)] $\displaystyle\sum_{k=1}^{\infty}\frac{1}{r_k}\E[ |\Delta M_k| : | \Delta M_k|  > \varepsilon r_k] < +\infty$ for any $\varepsilon > 0$, and
        \item [d)] $\displaystyle\sum_{k=1}^{\infty}\frac{1}{r_k^4}\E[ (\Delta M_k)^4 : |\Delta M_k|  > \delta r_k] < +\infty$ for some $\delta > 0$.
    \end{itemize}
    Then, $\displaystyle\limsup_{n\to\infty}\pm \frac{M_{\infty}-M_n}{\sqrt{2\cdot r^2_{n+1} \log |\log r^2_{n+1} |}} = 1$ a.s.
\end{itemize}
\end{lemma}

%
%%%%%%%%%%%%%%%%%%%%%%%%%%%%%%%%%%%%%%%%%%%%%%%%%%%%%%

%%%%%%%%%%%%%%%%%%%%%%%%%%%%%%%%%%%%%%%%%%%%%%%%%%%%%%%%%%%%%%%

\subsection{Proof of Theorem 1}
First, note that
\begin{equation}
    \frac{a_n}{n} = \frac{1}{n} \displaystyle\prod_{k=1}^{n-1}\left(1+\dfrac{\alpha}{k}\right)=\displaystyle\prod_{k=1}^{n-1}\left(\dfrac{k+\alpha}{k+1}\right)=\displaystyle\prod_{k=1}^{n-1}\left(1-\dfrac{1-\alpha}{k+1}\right).
\end{equation}
Moreover $-1\leq\alpha\leq1$, then  $0\leq\dfrac{k+\alpha}{k+1}\leq 1$. In addition, $\displaystyle\sum_{k=1}^{\infty}\frac{1-\alpha}{k+1}=\infty$ implies that $\lim \frac{a_n}{n}=0$. It follows that $(a_n)_{n\geq1}$ is non increasing. Let define $N_j=\frac{a_j}{j}\cdot\Delta M_j$. By \eqref{AcMn} we conclude that $(N_j)_{j\geq1}$ is a sequence of martingales such that
\begin{equation*}
    \displaystyle\sum_{j=1}^{\infty}\E[N_j^2\vert \mathcal{F}_{j-1}]\leq\displaystyle\sum_{j=1}^{\infty}\frac{2^2}{j^2}<\infty
\end{equation*}
From Theorem 2.17 in \cite{hall2014martingale}, $\displaystyle\sum_{j=1}^{\infty}N_j$ converges almost surely. Given that $\frac{n}{a_n}\rightarrow\infty$, by a direct application of Kronecker's lemma,
\begin{equation*}
    \frac{a_n}{n}\displaystyle\sum_{j=1}^{n}\Delta M_j=\frac{a_n}{n}M_n\overset{a.s}{\rightarrow}0
\end{equation*}
Then, we showed \eqref{lgn}.
On the other hand, using \eqref{lEsn}, and $\alpha<1$, we have that the first term is $o(n)$ as $n\rightarrow\infty$, and by \eqref{conv_an} and Lemma \ref{lemma_Sum_an}
\begin{equation*}
    \omega\cdot a_n \displaystyle\sum_{l=1}^{n-1}\frac{1}{a_{l+1}}\sim\omega\frac{n^\alpha}{\Gamma(\alpha+1)}\Gamma(\alpha+1)\frac{n^{1-\alpha}}{1-\alpha}=\frac{\omega }{1-\alpha}n
\end{equation*}
%Finalmente queda demostrado \eqref{limSn}
 
 %\begin{equation*}
 %\begin{array}{ll}
 %    \lim\limits_{n\to\infty}\frac{S_n}{n}=\frac{\omega}{1-\alpha} & a.s \\
 %\end{array}
 %\end{equation*}
%%%%%%%%%%%%%%%%%%%%%%%%%%%%%%%%%%%%%%%%%%%%%%%%%%%
\subsection{Proof of Theorem 2}

Let denote $\xi_n=X_n-\E(X_n\vert \mathcal{F}_{n-1})$, from \eqref{Var_M} we have $\E[\xi_{n}\vert\mathcal{F}_{n-1}]=0$ and
\begin{equation}
\label{E_epsi}
  \E[(\xi_{n+1})^2\vert\mathcal{F}_{n}]=\E[X_{n+1}^2\vert \mathcal{F}_{n}]-(\E(X_{n+1}\vert \mathcal{F}_{n}))^2,
\end{equation}
where
\begin{equation}
\label{EX2}
  \E[(X_{n+1})^2\vert\mathcal{F}_{n}]=\dfrac{n_+ + n_-}{n}\cdot (p+q)\cdot \theta+(1-\theta)\cdot (p+q).
\end{equation}
Let denote $Z_n=n_+ + n_-$, $\gamma=(p+q)\theta$ and $\tau=(1-\theta)\cdot(p+q)$. Then, using \eqref{XndadoF}
\begin{equation}
    \label{A4}
        \E[(\xi_{n+1})^2\vert\mathcal{F}_{n}]=\frac{\gamma}{n} Z_n+\tau-\left(\frac{\alpha}{n} S_n+\omega\right)^2.
\end{equation}
Similarly,
\begin{equation}
 \label{A5}
\begin{array}{l}
  \E[(\xi_{n+1})^4\vert\mathcal{F}_{n}]=\E\left[\left(X_{n+1}-\E(X_{n+1}\vert \mathcal{F}_{n})\right)^4\vert\mathcal{F}_{n}\right]\\[0.5cm]
      =\E\left[X^4_{n+1}-4X^3_{n+1}\E(X_{n+1}\vert \mathcal{F}_{n})+6X^2_{n+1}(\E(X_{n+1}\vert \mathcal{F}_{n}))^2\right.\\[0.5cm]
     \left.-4X_{n+1}\left(\E(X_{n+1}\vert \mathcal{F}_{n})\right)^3+\left(\E(X_{n+1}\vert \mathcal{F}_{n})\right)^4\vert\mathcal{F}_{n}\right]\\[0.5cm]
      =\E(X^2_{n+1}\vert\mathcal{F}_{n})-4(\E(X_{n+1}\vert\mathcal{F}_{n}))^2+6\E(X^2_{n+1}\vert\mathcal{F}_{n})\left(\E(X_{n+1}\vert\mathcal{F}_{n})\right)^2\\[0.5cm]
      -3\left(\E(X_{n+1}\vert \mathcal{F}_{n})\right)^4,
\end{array}   
\end{equation}

For any $n\geq1$, we have that $X_n$ is a random variable which takes values on $\{+1,0,-1\}$. In consequence, $\E[(\xi_{n+1})^2\vert\mathcal{F}_{n}]$ and $\E[(\xi_{n+1})^4\vert\mathcal{F}_{n}]$, are bounded and by \eqref{A4}, \eqref{A5} we have the almost sure limits
\begin{equation}
    \label{A6}
    \begin{array}{lll}
        \displaystyle\sup_{n\geq 0}~  \E[\xi^2_{n+1}\vert\mathcal{F}_{n}]\leq 2 & \text{ and } & \displaystyle\sup_{n\geq 0}~  \E[\xi^4_{n+1}\vert\mathcal{F}_{n}]\leq 16 \\
    \end{array}
\end{equation}
We deduce from \eqref{Var_M} and \eqref{A4} that
\begin{equation}
    \label{A7}
    \begin{array}{ll}
    \langle M\rangle_n & =\displaystyle\sum_{k=1}^{n}\frac{1}{a_k^2}\left[ \frac{\gamma}{k}Z_k+\tau-\left(\frac{\alpha}{k} S_k+\omega\right)^2\right]\\[0.5cm]
    \end{array}
\end{equation}
Now, note that using \eqref{PT1} and Theorem \ref{Teo1} we obtain $\E(X_{n+1}\vert \mathcal{F}_{n})  =\frac{\alpha}{n}S_n+\omega \rightarrow \alpha\frac{\omega}{1-\alpha}+\omega=\frac{\omega}{1-\alpha}$
Moreover, by \eqref{EX2} we have that $\E[X_{n+1}^2\vert \mathcal{F}_{n}]=\frac{\gamma}{n} Z_n+\tau$, in the Lemma \ref{l_mart} we constructed a martingale for $Z_n$, and in similar form than Theorem \ref{Teo1} it is possible to prove that:
\begin{equation}
\label{limZn}
    \lim\limits_{n\to\infty}\frac{Z_n-\E(Z_n)}{n}=0 \text{ and } \lim\limits_{n\to\infty}\frac{\E(Z_n)}{n}=\frac{\tau}{1-\gamma}
\end{equation}
Then,
\begin{equation}
    \label{convX2}
    \E[X_{n+1}^2\vert \mathcal{F}_{n}]  = \frac{\gamma}{n} Z_n+\tau \sim \frac{\tau}{1-\gamma}
\end{equation}
Finally, from \eqref{E_epsi}
\begin{equation}
\label{var}
   \E[(\xi_{n+1})^2\vert\mathcal{F}_{n}]  \sim \frac{\tau}{1-\gamma}-\left(\frac{\omega}{1-\alpha}\right)^2= \phi
\end{equation}
as defined in \eqref{notation}. Given that $\Delta M_k = \frac{\xi_k}{a_k}$, the expectation of the square of the martingale difference satisfies,
\begin{equation}
\label{varM}
\begin{array}{ll}
   \E[(\Delta M_k)^2\vert\mathcal{F}_{k-1}] &\sim \frac{\phi}{a_k^2}
\end{array}
\end{equation}

\begin{center}
    \textbf{Diffusive region, $0\leq\alpha<1/2$}
\end{center}

We use central limit theorem for martingales as stated in Corollary 2.1.10 in \cite{duflo2013random}. First, note that from \eqref{vn} and using Lemma \ref{lan2}, we have that $v_n$ diverges in this region. Moreover, by \eqref{Var_M} and \eqref{varM} we obtain that
\begin{equation}
    \label{lim_varcua/vn}
    \begin{array}{cc}
      \lim\limits_{n\to\infty} \frac{\langle M\rangle_n}{v_n}=\phi    &  a.s\\
    \end{array}
\end{equation}
Now, we need to prove that $(M_n)$ satisfies the Lindeberg's condition, that is, for any $\varepsilon>0$,
 \begin{equation}
     \label{con_lindebe}
     \frac{1}{v_n}\displaystyle\sum_{k=1}^{n}\E\left[\vert\Delta M_{k+1}\vert^2 \mathbb{I}_{\vert\Delta M_{k+1}\vert\geq\varepsilon\sqrt{v_n}}\vert\mathcal{F}_n\right]\overset{P}{\rightarrow}0,
 \end{equation}
 it can be seen, for instance, in Theorem 2.1.9 from \cite{duflo2013random}. Then using \eqref{A6}, we have:
 \begin{equation*}
     \begin{array}{l}
    \displaystyle\frac{1}{v_n}\displaystyle\sum_{k=1}^{n}\E\left[\vert\Delta M_{k+1}\vert^2 \mathbb{I}_{\vert\Delta M_{k+1}\vert\geq\varepsilon\sqrt{v_n}}\vert\mathcal{F}_n\right]       \leq \frac{1}{\varepsilon^2v_n^2}\displaystyle\sum_{k=1}^{n}\E\left[\vert\Delta M_{k+1}\vert^4 \vert\mathcal{F}_n\right] \\[0.5cm]
           \leq \displaystyle\sup_{1\leq k\leq n}\E\left[\xi^4_{k+1} \vert\mathcal{F}_n\right]\frac{1}{\varepsilon^2v_n^2}\displaystyle\sum_{k=1}^{n} \frac{1}{a_k^4} \leq\frac{16}{\varepsilon^2v_n^2}\displaystyle\sum_{k=1}^{n} \frac{1}{a_k^4}\\[0.5cm]
     \end{array}
 \end{equation*}

Now, note that $a_n^{-4}v_n^{-2}$ is equivalent to $(1-2\alpha)n^{-2}$ and it is well known that $\sum_{n=1}^{\infty} \frac{1}{n^2}=\frac{\pi^2}{6}$. Therefore
\begin{equation}
    \label{a_4/vn}
   \displaystyle\sum_{n=1}^{\infty} \frac{1/a_n^4}{v_n^2}<\infty
\end{equation}
By the Kronecker's lemma (see \cite{duflo2013random}), we obtain
\begin{equation}
    \label{lim_vn_ak4}
    \lim\limits_{n\to\infty}\frac{1}{v_n^2}\displaystyle\sum_{k=1}^{n} \frac{1}{a_k^4}=0
\end{equation}
In this sense, by central limit theorem for martingales, we have
\begin{equation}
\label{mclt}
    \lim\limits_{n\to\infty}\frac{M_n}{\sqrt{v_n}}\overset{d}{\rightarrow}N(0,\phi)
\end{equation}
Since $M_n=\frac{S_n-\E(S_n)}{a_n}$, by \eqref{an} and \eqref{conv_dif_vn}, the equation \eqref{TCLdif} holds.
Finally by \eqref{lim_varcua/vn} and given that $v_n/v_{n+1} \to 1$, as $n \to\infty$, we use Lemma 3.5 from \cite{JJQ} to obtain,
\begin{equation}
\label{LILproof}
\displaystyle\limsup_{n\to\infty} \pm \frac{M_n}{\sqrt{2\phi v_n\log \log (\phi v_n)}} = 1 \ \text{ a.s. }
\end{equation}

\begin{center}
    \textbf{Critical region, $\alpha=1/2$}
\end{center}
In this region we proceed in a similar way than previous case. Note that $v_n$ diverges and using Lemma \ref{lan2}, and \eqref{A7} we obtain \eqref{lim_varcua/vn}.
Moreover \eqref{a_4/vn} holds since $a_n^{-4}v_n^{-2}$ is equivalent to $(n \log n)^{-2}$ by \eqref{conv_crit_vn} and it is well known that $\sum_{n=1}^{\infty} \frac{1}{(n\log n)^2}<\infty$.
We realized in the diffusive region that \eqref{a_4/vn} is a sufficient condition to have $(M_n)$ satisfying the Lindeberg's condition. In consequence, we use the central limit theorem for martingales \eqref{mclt}, equation \eqref{Mart} and \eqref{conv_crit_vn} to obtain
\begin{equation*}
    \lim\limits_{n\to\infty}\frac{M_n}{\sqrt{v_n}}=\frac{S_n-\E(S_n)}{\sqrt{n \log n}}
\end{equation*}
Then \eqref{TCLcrit} holds. Moreover, to prove \eqref{LIL2} we use Lemma 3.5 from \cite{JJQ}, we obtain the analogous of \eqref{LILproof}, given \eqref{conv_crit_vn} we complete the proof.

%\begin{equation*}
 %   \frac{s_n-\E(S_n)}{\sqrt{\phi\cdot n \log n}}\overset{d}{\rightarrow}N(0,1)
%\end{equation*}

 %%%%%%%%%%%%%%%%%%%%%%%%%%%%%%%%%%%%%%%%%%%%%%%%%%%%%%%%%%%%%%%
 \subsection{Proof of Theorem 3}
 In superdiffusive region, $1/2<\alpha\leq1$, we use Lemma \ref{Kubotalem}. Note that by \eqref{varM} and the bounded convergence theorem, $\displaystyle\sum_{k=1}^{\infty}\E[(\Delta M_k)^2]\sim \phi\Gamma(\alpha+1)^2\displaystyle\sum_{k=1}^{\infty}\frac{1}{k^{2\alpha}}$. Then, if $\alpha> 1/2$, we have $\displaystyle\sum_{k=1}^{\infty}\E[(\Delta M_k)^2] < \infty$. Therefore, from item (i) of the cited Lemma, the following limit exists with probability one,
\begin{equation*}
       W:=\displaystyle\sum_{k=1}^{\infty}\Delta M_k =\lim\limits_{n\to\infty}  \frac{S_n-\E[S_n]}{a_n},
\end{equation*}
and $M_n \overset{L^2}{\rightarrow} W$. Then $\E(W)=0$ and $\E(W^2)= \sum_{k=1}^{\infty}\E[(\Delta M_k)^2] > 0$.\\

To prove \eqref{T_cnv_w} we use item ii) of Lemma \ref{Kubotalem}, for that, we need to prove for $\alpha>1/2$ that items a) and b) hold.\\
A sufficient condition for a) is:
   \begin{equation}
     \label{Cond_suf_a)}
     \frac{1}{r_n^2}\sum_{k=n}^{\infty}\E\left[(\Delta M_k)^2\vert\mathcal{F}_{k-1}\right]\rightarrow 1~\text{ as }~ n\to\infty ~\text{ in probability}
 \end{equation}
by Corollary 1 item b`) in \cite{heyde1977central}, where $r_n^2=\sum_{k=n}^{\infty}\E[(\Delta M_k)^2]$ and given \eqref{varM}, we have:

\begin{equation*}
\begin{array}{ll}
     \displaystyle\sum_{k=n}^{\infty}\E[(\Delta M_k)^2\vert\mathcal{F}_{k-1}]&=\displaystyle\sum_{k=n}^{\infty}\frac{\phi}{a_k^2}
    \sim \phi\Gamma(\alpha+1)^2\displaystyle\sum_{k=n}^{\infty}\frac{1}{k^{2\alpha}}\\[0.4cm]
    &\sim\displaystyle\frac{\phi\Gamma(\alpha+1)^2}{(2\alpha-1)n^{2\alpha-1}}
    \sim \frac{\phi n}{(2\alpha-1)a_n^2}~a.s
\end{array}
\end{equation*}
Using the bounded convergence theorem, we obtain that:
\begin{equation*}
     \sum_{k=n}^{\infty}\E[(\Delta M_k)^2]=\sum_{k=n}^{\infty}\frac{\phi}{a_k^2}\sim \frac{\phi n}{(2\alpha-1)a_n^2}~a.s
\end{equation*}
Then, as $n\to\infty$, we have:
\begin{equation*}
     \displaystyle\lim_{n\to\infty}\frac{1}{r_n^2}\sum_{k=n}^{\infty}\E\left[(\Delta M_k)^2\vert\mathcal{F}_{k-1}\right]=1~a.s
 \end{equation*}
This shows that conditions a) and a’) in Lemma \ref{Kubotalem} are satisfied. Moreover
 \begin{equation*}
     r_n^4\sim \frac{\Gamma(\alpha+1)^4\phi^2 }{(2\alpha-1)^2n^{4\alpha-2}} \text{ and } (\Delta M_n)^4\leq\frac{16}{(a_n^4)}\sim \frac{16\Gamma(\alpha+1)^4}{n^{4\alpha}}.
\end{equation*}
Now, the sufficient condition to have b) is the Remark A1 in \cite{kubota2019gaussian}, in this sense, using Theorem 2.1.9 from \cite{duflo2013random} and \eqref{A6}:
\begin{equation*}
\begin{array}{ll}
         \frac{1}{r_n^2}\displaystyle\sum_{k=n}^{\infty}\E\left[(\Delta M_{k+1})^2: \vert\Delta M_{k+1}\vert\geq\varepsilon r_k\right]  & \leq \frac{1}{\varepsilon^2r_n^4}\displaystyle\sum_{k=n}^{\infty}\E\left[\vert\Delta M_{k+1}\vert^4\right] \\[0,4cm]
         &\leq \displaystyle\frac{16\Gamma(\alpha+1)^4}{\varepsilon^2r_n^4}\sum_{k=n}^{\infty} \frac{1}{a_k^4}
         \leq\frac{C_1}{r_n^4}\displaystyle\sum_{k=n}^{\infty} \frac{1}{a_k^4}
\end{array}
 \end{equation*}
At this point, we need to prove that
\begin{equation}
    \label{lim_rn_ak4}
    \lim\limits_{n\to\infty}\frac{1}{r_n^4}\displaystyle\sum_{k=n}^{\infty} \frac{1}{a_k^4}=0
\end{equation}
which clearly holds since $a_n^{-4}r_n^{-4}$ is equivalent to $\frac{1}{a_n^2\phi}$ and by Lemma \ref{lan2}, in the superdiffusive case $\displaystyle\sum_{n=1}^{\infty} \frac{1}{a_n^2}<\infty$. Therefore, $\displaystyle\sum_{n=1}^{\infty} \frac{1/a_n^4}{r_n^4}<\infty$
and by Lemma 1 item ii) in \cite{heyde1977central}, \eqref{lim_rn_ak4} holds. Then, the proof finishes by using the Theorem 1 (b) in \cite{heyde1977central} 
\begin{equation*}
    \frac{W-M_n}{r_{n+1}}\overset{d}{\rightarrow}N(0,1)
\end{equation*}
Finally, for $\varepsilon>0$ we have that:
\begin{equation*}
         \frac{1}{r_k}\E\left[\vert\Delta M_{k+1}\vert: \vert\Delta M_{k+1}\vert\geq\varepsilon r_k\right] \leq\frac{1}{r_k} \frac{1}{\varepsilon^3r_n^3}\E\left[\vert\Delta M_{k+1}\vert^4\right]
         \leq C_3 k^{4\alpha-2}k^{4\alpha}
         \leq\frac{C_3}{k^2}
 \end{equation*}
which implies that c) holds. In the case of d), it is enough to note that $\sum_{n=1}^{\infty} \frac{1}{r_k^4}\E[(\Delta M_k)^4]<\infty$. Then, \eqref{LIL3} holds.

%%%%%%%%%%%%%%%%%%%%%%%%%%%%
%%%%%%%%%%%%%%%%%%%%%%%%%%%%
%%%%%%%%%%%%%%%%%%%%%%%%%%%%

\section*{Acknowledgements}
The authors thank the anonymous referees for their valuable comments. Also several corrections by Christian Caama\~no and Rodrigo Lambert. This work was partially supported by Fondecyt nº 11200500

\medskip

{\scriptsize
{\sc Departamento de Estad\'i{}stica, Universidad del B\'io-B\'io. Avda. Collao 1202, CP 4051381, Concepci\'on, Chile. E-mail address: magonzalez@ubiobio.cl
}
}

\end{document}